%% file: convexQLearning_09_14_EXTENDED_spm.tex
\documentclass[11pt]{article}   % Comment this line out if

\usepackage{geometry}
\pdfoutput=1  %  for arxiv

\usepackage[usenames,dvipsnames]{color}
\usepackage{color}

\usepackage{amsmath, amssymb, graphicx, url, algorithm2e}

\usepackage[colorlinks,%
linkcolor=BrickRed,% 
filecolor =BrickRed,%
urlcolor=blue,% 
citecolor=RoyalPurple,%
]{hyperref}

\usepackage{flushend}

\usepackage{optidef}
\usepackage{graphics} % for pdf, bitmapped graphics files
\graphicspath{{figures/}}

\usepackage{url}            % simple URL typesetting

\usepackage{booktabs}       % professional-quality tables

\usepackage{wrapfig}

\usepackage{caption}
\usepackage{marginnote}
\usepackage{float} 

\usepackage{algpseudocode}

\usepackage{textgreek,upgreek,bm}

\usepackage{ulem} 

\input{bookmacrosRL_IEEE}

\def\wham#1{\smallbreak\pagebreak[3]%
	\noindent\textbf{#1}\ \ \gobblepars}

\def\avgpsi{\bar{\psi}^{\upmu}}

\def\Treset{T} 
\def\sTdiff{\widehat{\mathcal{D}}} % state dependent sampling TD
\def\sTdiffDQN{\widehat{\mathcal{D}}^{\text{DQN}}} % state dependent sampling TD

\def\nU{n_U}
\def\clA{\mathcal{A}}
\def\clC{\mathcal{C}}

\def\regBCQ{\mathcal{G}}

\def\varC{\mathcal{C}}

\graphicspath{{figures/}}

\newlength{\noteWidth}
\long\def\notes#1{\ifinner
	{\tiny #1}
	\else
	\marginpar{\parbox[t]{\noteWidth}{\raggedright\tiny #1}}
	\fi}

        \def\notes#1{\typeout{See notes!}}

\def\sfb#1{}

\newtheorem{theorem}{Theorem}[section]

\newtheorem{proposition}[theorem]{Proposition}
\newtheorem{lemma}[theorem]{Lemma}

\usepackage{cleveref}

\Crefname{corollary}{Corollary}{Corollaries}
\Crefname{eqnarray}{eq.}{eqs.}
\Crefname{equation}{eq.}{eqs.}

\Crefname{figure}{Fig.}{Figs.}
\Crefname{tabular}{Tab.}{Tabs.}
\Crefname{table}{Tab.}{Tabs.}
\Crefname{proposition}{Prop.}{Propositions}
\Crefname{theorem}{Thm.}{Thms.}
\Crefname{definition}{Def.}{Defs.} 
\Crefname{section}{Section}{Sections}
\Crefname{lemma}{Lemma}{Lemmas}
\Crefname{assumption}{Assumption}{Assumptions} 

\def\psisub#1{\psi_{(#1)}}

\def\csub#1{c_{(#1)}}

\def\ctheta{\check{\theta}}

\makeatletter
\def\thanks#1{\protected@xdef\@thanks{\@thanks
		\protect\footnotetext{#1}}}
\makeatother
\title{\LARGE \bf
	Sufficient Exploration for Convex Q-learning}

\author{Fan Lu, Prashant G. Mehta, Sean P. Meyn  and Gergely Neu%<-this % stops a space
	\thanks{SPM and FL are with the Department of Electrical and Computer Engineering, University of Florida, Gainesville, FL 32611;   SPM holds an Inria International Chair, Paris, France.  
	}%
	\thanks{PGM is with the Coordinated Science Laboratory and the Department of Mechanical Science and Engineering at the University of Illinois at Urbana-Champaign (UIUC).}
	\thanks{GN is with the Department of Information and Communication Technologies, Universitat Pompeu Fabra (Barcelona, Spain).}
	\thanks{ Financial support from ARO award W911NF2010055 and NSF
		award  EPCN 1935389}    
}

\begin{document}
	
	\maketitle

	\begin{abstract}%
		In recent years there has been a collective research effort to find new formulations of  reinforcement learning that are simultaneously more efficient and more amenable to analysis.  This paper concerns one approach that builds on the linear programming (LP) formulation of optimal control of Manne.    A primal version is called logistic Q-learning, and a dual variant is convex Q-learning.   This paper focuses on the latter, while building bridges with the former.  The main contributions follow: 
		(i) The dual of convex Q-learning is not precisely  Manne's LP or a version of logistic Q-learning, but has similar structure that reveals the need for regularization to  avoid over-fitting.   
		(ii) A sufficient condition is obtained for a bounded solution to the Q-learning LP.   
		(iii)   Simulation studies reveal numerical challenges when addressing sampled-data systems based on a continuous time model.  The challenge is addressed using state-dependent sampling.   The theory is illustrated with applications to examples from OpenAI gym. It is shown that convex Q-learning is successful in cases where standard Q-learning diverges, such as the LQR problem.

		 {\color{blue}
		 Note:  this is a slightly extended version of CDC 2022 \cite{lumehmeyneu22}, including an Appendix with proofs of the main results.}
	\end{abstract}

	\section{Introduction}
	\label{s:intro}

	Ever since the introduction of Watkins' Q-learning algorithm in the 1980s, the research community has searched for a general theory beyond the so-called tabular settings (in which the function class spans all possible functions of state and action).    The natural extension of Q-learning to  general  function approximation setting seeks to solve what is known as a \textit{projected Bellman equation} (PBE).   There are few results available giving sufficient conditions for the existence of a solution,  or convergence of the algorithm if a solution does exist  \cite{sutszemae08,melmeyrib08,leehe19}.
	Counterexamples show that conditions on the function class are required in general,  even in a linear function approximation setting  \cite{bai95,tsivan96,gor00}.   The GQ-algorithm of
	\cite{maeszebhasut10} is one success story,  based on a relaxation of the PBE.
	
%	\notes{is sut95 valuable?}
	
	Even if existence and stability of the algorithm were settled, we would still face the challenge of interpreting the output of a Q-learning algorithm based on the PBE criterion.   Inverse dynamic programming provides bounds on performance, but only subject to a weighted $L_\infty$ bound on the Bellman error, while RL theory is largely based on $L_2$ bounds \cite{sutbar18,CSRL}.
	
	Both logistic Q-learning \cite{bascurkraneu21}  and convex Q-learning \cite{mehmey09a,leehe19b,mehmeyneulu21} are based on the convex analytic approach to optimal control of \cite{man60a} and its significant development over the past 50 years in the control and operations research literature \cite{bor02a,schsei85,wanboy09,farroy06}.   
	There is a wealth of unanswered questions:  
	
	%\archive{Perhaps most important, but \textit{not} covered in this paper is an explanation for the empirical success of Q-learning.   It seems that stability is typical, even without theory  (but we'll keep this to ourselves}
	
	\wham{(i)}
	It is shown in a tabular setting that the dual of convex Q-learning is somewhat similar to Manne's primal LP~\cite{man60a}, but its sample path form also brings differences. 
	
	\wham{(ii)}
	The most basic version of convex Q-learning is a linear program (LP).  It is always feasible, but boundedness has been an open topic for research (except for stylized special cases).    It is shown here for the first time that boundedness holds if the covariance associated with the basis is full rank.    This may sound familiar to those acquainted with the literature, but the proof is non-trivial since  theory is far from the typical $L_2$ setting of TD-learning  \cite{tsivan96,sutbar18,CSRL,bertsi96a}.   
	
	%	\archive{we did not establish equivalence.   And remember,   boundedness of the constraint region is not necessary for boundedness of the solution to the LP! }
	
	So far,   LP formulations of reinforcement learning (RL) have restricted to either the tabular or `linear MDP'' settings \cite{bascurkraneu21,leehe19b},   or deterministic optimal control with general linear function approximation \cite{mehmeyneulu21}.   In this paper we focus on the latter, since the challenges in the stochastic setting would add significant complexity.
	
	We consider a nonlinear state space model in discrete time,
	\begin{equation}
		x(k+1) = \Dds(x(k),u(k)) \,,\qquad k\ge 0\,, \  x(0) =x \in\state\, .
		\label{e:dis_NLSSx_controlled_opt}
	\end{equation}
	The state space  $\state$ is a closed subset of $\Re^n$, and the input (or action) space $\ustate$ is finite,   with cardinality $\nU \eqdef |\ustate|$,  and where $\Dds\colon  \state\times\ustate \to\state $.     We may have state-dependent constraints, so that for each $x\in\state $ there is a set $\ustate(x)\subset \ustate$ for which $u(k)$ is constrained to $\ustate(x(k))$ for each $k$.
	Notation is simplified by denoting  $\{ z(k) = (x(k),u(k)) : k\ge 0\}$, evolving on  $\zstate \eqdef \{(x,u) : x\in\state\,, \,  u\in\ustate(x) \}$.   
	
	\begin{figure*}
		\centering
		\includegraphics[width=0.95\hsize]{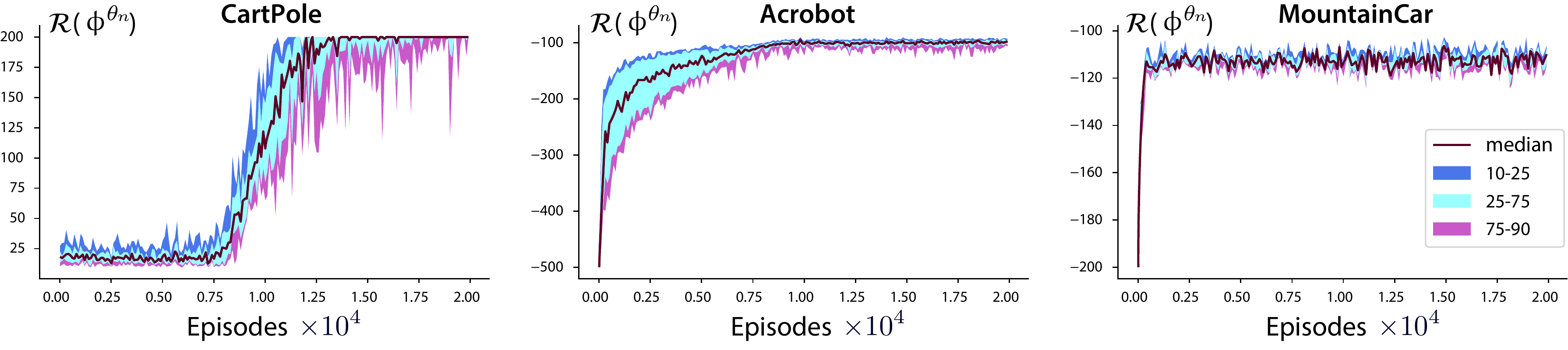}
		\caption{Average cumulative reward as a function of iteration for three examples.   Percentiles are estimated using independent runs.}
		\label[figure]{f:allresults}
	\end{figure*}	
	
	The paper concerns infinite-horizon optimal control, 
	whose definition requires a cost function $c\colon \zstate\to\Re_+$,  and a pair $z^e \eqdef (x^e,u^e)\in\zstate$ that achieves equilibrium:
	\[
	x^e = \Dds (x^e,u^e).
	\] 
	The cost function $c\colon\zstate\to\Re_+$ vanishes   at $z^e$.

	These assumptions are imposed so that there is hope that the (optimal) Q-function is finite valued:   
	\begin{equation}
		\begin{aligned}
			Q^\oc (x,u) 
			&= \min_{u(1),u(2),\dots} \sum_{k=0}^\infty  c (x(k) , u(k))    \,, 
			\\
			x(0)& =x\in\state\,,\  u(0) = u \in\ustate (x)    
		\end{aligned}
		\label{e:Q_part1}
	\end{equation} 
%	The goal of optimal control is to find an optimizing input sequence.   %, and in the process we often need to compute the Q-function.	
	The Bellman equation may be expressed
	\begin{equation}
		Q^\oc(x,u)  =  c (x,u) 
		+   \uQ^\oc ( \Dds(x,u) ),
		\label{e:WeAreQ}
	\end{equation}
	with $ \uQ(x) = \min_u Q(x,u)$ for any function $Q$.  The optimal input is state feedback $u^\oc(k) = \fee^\oc(x^\oc(k))$, using the ``$Q^\oc$-greedy'' policy, 
	\begin{equation}
		\fee^\oc (x) \in \argmin_{u \in \ustate(x)}  Q^\oc(x,u)\,,\qquad x\in\state.
		\label{e:fee}
	\end{equation}
	
	Q-learning algorithms are designed to approximate $Q^\oc$ within a parameterized family of functions $\{ Q^\theta : \theta\in\Re^d\}$,  and based on an appropriate approximation obtain a policy in analogy with \eqref{e:fee}: 	
	\begin{equation}
		\fee^\theta(x) \in \argmin_{u \in \ustate(x)} Q^\theta(x,u)
		\label{e:Qgreedy}
	\end{equation}

For any input-state sequence $\{ u(k), x(k) : k\ge 0\}$, the Bellman equation \eqref{e:WeAreQ} implies
\begin{equation}
		Q^\oc(z(k) )  =  c (z(k))   +   \uQ^\oc ( x(k+1))
		\label{e:e:WeAreQ2}
\end{equation}
This motivates the 	\textit{temporal difference sequence}:
for any $\theta$,    the observed error at time $k$ is denoted
	\begin{equation}
		\Tdiff^\circ_{k+1}(\theta) \eqdef   -  Q^\theta (z(k) )   +   c(z(k))  +   \uQ^\theta (x(k+1)) 
		\label{e:TDonline}
	\end{equation} 
	Given observations over a time-horizon $0\le k\le\Nsam$, one approach is to choose $\theta^\ocp$ that minimizes the  \textit{mean-square Bellman error}:  
	\begin{equation*}
		\frac{1}{\Nsam} \sum_{k=0}^{\Nsam-1}  \bigl[   \Tdiff^\circ_{k+1}(\theta)  \bigr]^2
		%\label{e:MSBE} 
	\end{equation*}
	The GQ-algorithm is designed to solve a similar non-convex optimization problem.
	There has been great success using an alternative \textit{Galerkin relaxation} \cite{CSRL}:  	A sequence of $d$-dimensional \textit{eligibility vectors} $\{ \elig_k  \} $ is constructed,  and the goal then is to solve a version of the projected Bellman equation,
	\begin{equation}
		\Zero  =  \frac{1}{\Nsam} \sum_{k=0}^{\Nsam-1}\Tdiff^\circ_{k+1}(\theta) \elig_k 
\label{e:EmpiricalGalerkin}
\end{equation} 
%\rd{
%CHANGE N TO K FOR CONSISTENCY
%\\
% $\Tdiff_{n+1}  =  \Tdiff_{n+1} ^\circ(\theta_n)
%$
%\begin{equation}
%		\theta_{n+1} = \theta_n + \alpha_{n+1} \Tdiff_{n+1}  		\psi(x_n, u_n).
%\label{e:Q}
%\end{equation}}

The  standard Q-learning algorithm  is based on the temporal difference sequence $\Tdiff_{k+1} \eqdef  \Tdiff_{k+1} ^\circ(\theta_k)$ in the recursion 
\begin{equation}
	\theta_{k+1} = \theta_k + \alpha_{k+1} \Tdiff_{k+1} \elig_k^{\theta_k}
	\label{e:Q}
\end{equation}
  with $\elig_k^{\theta_k} = \nabla_\theta Q^\theta(z(k))|_{\theta = \theta_k}$,  and $\{\alpha_{k+1}\}$ is a non-negative step-size sequence
 \cite{sutbar18,CSRL}.  When convergent, the limit satisfies a  projected Bellman equation similar to \eqref{e:EmpiricalGalerkin}:
	\begin{equation}
		\begin{aligned}
			\barf(\theta^\ocp) &= 0 \,,\qquad 
			\barf(\theta) =
			\Expect
			\bigl[      \Tdiff^\circ_{k+1}(\theta) \elig_k^{\theta}    \bigr]   
		\end{aligned} 
		\label{e:WatkinsRelax}
\end{equation}
While not obvious from its description, parameter estimates obtained from the DQN algorithm solve the same
projected Bellman equation, provided it is convergent  (see \cite{mehmeyneulu21}).

%\notes{OK for paper, but no need here: (a significant component of famous applications such as AlphaGo, cite{silver2016mastering})   }
%\notes{
%FL2SM: projected bellman equation? should we cite some equation? Is it the right handside of \eqref{e:WatkinsRelax}.
%\\
%Is this resolved?	}
		
	\smallbreak

	%The starting point of this paper is to clarify the paradox that approximating the Q-function involves solving a nonlinear fixed point equation \eqref{e:WeAreQ} or \eqref{e:e:WeAreQ2},  for which root finding problems for approximation may not be successful, and from this create a new family of RL algorithms designed to minimize a convex variant of the empirical mean-square error \eqref{e:MSBE}.  This step was anticipated in   \cite{mehmey09a}, for which a convex formulation of Q-learning was proposed for deterministic systems in continuous time. 
	
	The main contributions and organization of the paper are summarized as follows:
	
	%archive{SM2FL, note my macro `wham'.   It is much better than the previous 'head' because there is no risk of a paragraph break.}
	
	\wham{(i)} 
	
	In TD-learning it is known that the basis must be linearly independent to obtain a unique solution.   Theory developed in \Cref{s:FanTheory} implies that a similar condition is both necessary and sufficient to obtain a bounded constraint region in the convex program that defines convex Q-learning.   This result is obtained in the general setting with linear function approximation, so in particular the state space need not be finite.   The main conclusions are summarized in 
	\Cref{t:FanTheorem2}.

	\begin{wrapfigure}[14]{r}{0.4\textwidth}
 		\begin{center}
			\includegraphics[width=.95\hsize]{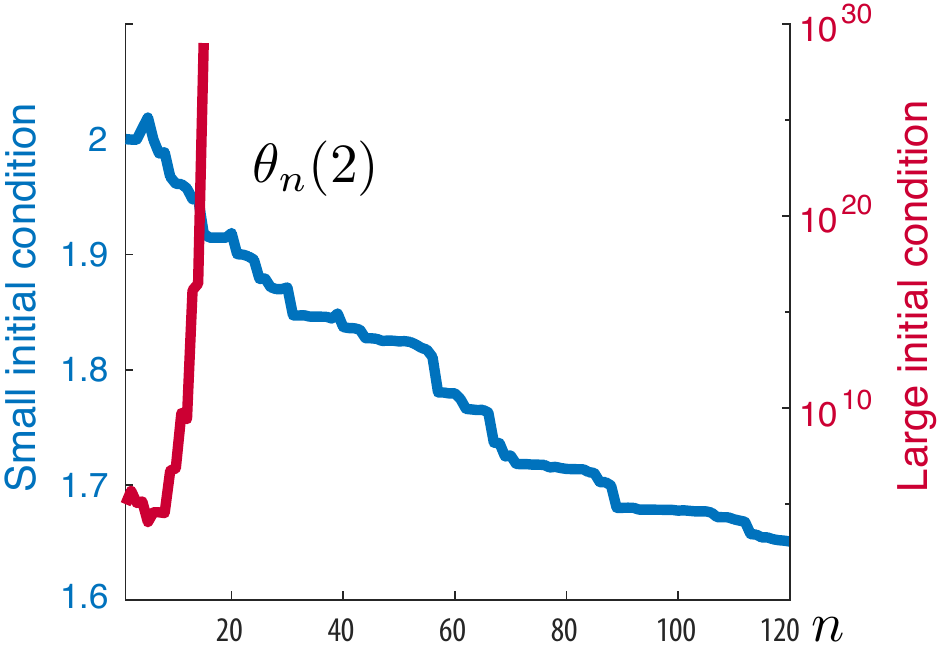}
		\end{center}
		\vspace{-1em}
		
		\caption{Evolution of one parameter using the standard Q-learning algorithm  \eqref{e:Q} in an application to the linear quadratic regulator problem.
			 }
		\label[figure]{f:compCvxQDQN}
	\end{wrapfigure}
	
	\notes{I deleted this "The estimates converge only when the norm of the initial condition is sufficiently small".    I believe it is false -- the estimates never converge!}

	\notes{sm2fl:   note use of hsize and not textwidth.  
	Also, note that I had a wrapfig version of the figure.}
	
	\wham{(ii)}
	The dual of convex Q-learning is described in  \Cref{sec:dual}, along with a number of consequences:  \Cref{t:FanPrimalDual} provides an interpretation of complementary slackness as an exact solution to the dynamic programming equation at selected state-action pairs.  This suggests that regularization is needed to avoid over-fitting in general.
	
	In the tabular setting,   the rank condition ensuring a bounded constraint region is equivalent to full exploration of all state-input pairs;   this is also a sufficient condition to ensure that convex Q-learning will compute exactly the optimal Q-function.  In this special case, the dual is similar to the primal introduced by Manne \cite{man60a}.	
	%archive{SM2FL: did you see this? FL: .
		%\\
		%There were errors in previous statement, copied here:
		%\\	In the tabular setting,   a bounded primal requires full exploration and hence convex Q-learning will compute exactly the optimal Q-function.
		%\\
		%	First, a technical error.   It is the constraint region that is bounded provided we have a full rank condition.  The word ``requires'' is just incorrect.    Second,  there are two concepts in a single sentence!  That is, boundedness and exact computation
		%\\
		%	Also, I deleted this since it is not easy to justify here:  , and also to the primal objective in    logistic Q-learning.
		%}

	\wham{(iii)}   Simulation studies reveal numerical challenges when addressing sampled-data systems;   the challenge is addressed here using state-dependent sampling.

%	\begin{figure}[h]
%		\centering
%		\includegraphics[width=0.35\hsize]{figures/compCvxQDQN_1DLQR_flat} 
%		\caption{Evolution of one parameter using standard-Q learning.
%		The estimates converge only when the norm of the initial condition is sufficiently small.   }
%		\label[figure]{f:compCvxQDQN}
%	\end{figure}

	\wham{(iv)}  Theory is illustrated in  \Cref{s:nres} with applications to examples from OpenAI gym. 
	\Cref{f:allresults} shows average cumulative reward as a function of training time for three examples.
		 The algorithm is remarkably robust, and is successful in cases where standard Q-learning diverges, such as the LQR problem.   An example of divergence is shown  in \Cref{f:compCvxQDQN}.    Details may be found in 	\Cref{s:nres}.

	% SM2FL:  This is an odd place for this statement:
	% The ODE of DQN is consistent with that of standard Q-learning \cite{CSRL, sutbar18}.}.  

%	SM2FL:	I don't understand what is asked of me here:
%	\choreL{We need to add the comparison of cvxQ with DQN figure and refer to it here. }
%   For journal we need real DQN
	
	%\rd{This looks really interesting -- merge this into the earlier pages:}
	%The concurrent work \cite{bascurkraneu21} bears several similarities to the present paper: 
	%most notably, they present another convex alternative to the squared Bellman error~\eqref{e:MSBE} called the 
	%logistic Bellman error, which is also derived from the same DPLP we use as our starting point. That said, the two papers 
	%diverge significantly when it comes to the details of the algorithm design and the nature of the performance 
	%guarantees. Unifying the advantages of the two approaches is an exciting direction for future research.
	%\choreL{\textbf{TODO:} Paragraph on organization of the paper. This can be done after the paper is finalized.}
	%\textbf{The remainder of the paper is organized as follows:} 

	\section{Convex Q-learning}
	\label{s:OPTSA} 
	
	Convex Q-learning is motivated by the classical LP characterization of the value function.    
	For any   $Q\colon\zstate\to\Re$, and     $r\in\Re$,  let $S_Q(r)$ denote the \textit{sub-level set}:
	\begin{equation*}
		S_Q(r) =  \{ z \in\zstate : Q(z) \le r\}
		\label{e:SJ}
	\end{equation*}
	The function $Q$ is called \textit{inf-compact} if the set $S_Q(r)$ is   pre-compact for $r$ in the range of $Q$.   The following may be found in \cite[Ch.~4]{CSRL}.
	
	%FL2SM: I REMOVE $z^e \eqdef (x^e, u^e)$ since $z^e$ has already been defined in the introduction. 
	
	\begin{subequations}
		\begin{proposition}
			\label[proposition]{t:LPforJQ}
			Suppose that the value function $Q^\oc$ defined in \eqref{e:Q_part1} is continuous, inf-compact, and vanishes only at $z^e$.   Then, for any positive measure $\upmu$ on $\state\times\ustate$,   $Q^\oc$  solves  the following convex program:  
			\begin{align} 
				\!\!\!			\max_{Q}   \ \  &   \langle \upmu, Q \rangle      
				\label{e:ConvexBEa}
				\\
				\st  \ \  &  Q(z) \le  c (z)  +  \uQ( \Dds(z) )  \,, \quad  z\in\zstate 
				\\
				\ \ & \text{$Q$ is continuous, and $ Q(z^e) =0$. }
				\label{e:ConvexBEb} 
			\end{align} 
		\end{proposition}%
		The nonlinear operation that defines $\uQ$ can be removed if the input space $\ustate$ is finite, so that \eqref{e:ConvexBE} can be represented as an LP;  it is always a convex program (even if infinite dimensional) because this minimization operator is a concave functional.  The LP construction is based on the inequalities in \eqref{e:newconstt} below.
		
			\label{e:ConvexBE}%
	\end{subequations}%

	Convex Q-learning is based on an approximation of the  convex program \eqref{e:ConvexBE},  
	seeking an approximation to $Q^\oc$ among a finite-dimensional family  $\{ Q^\theta : \theta\in\Re^d\}$.  
	The value $\theta_i$ might represent the $i$th weight in a neural network function approximation architecture,  but to justify the adjective \textit{convex} we require a linear family:
	\begin{equation}
		Q^\theta(z) = \theta^\transpose \psi (z) 
		\label{e:linParFamily}
	\end{equation}
	subject to the constraint $\psi_i(z^e) = 0$, for each $1\le i \le d$. 
	On introducing the $d$-dimensional vector
	\begin{equation}
		\avgpsi \eqdef \sum_{ z\in\zstate} \upmu(z) \psi(z) = \langle \upmu, \psi \rangle
		\label{e:mupsi}
	\end{equation}	
	it follows that $ \langle \upmu, Q^\theta \rangle       = \theta^\transpose \avgpsi$.
	
	Consider the restriction of \eqref{e:ConvexBE} to this parameterized family:
	\begin{subequations}
		\begin{align} 
			\!\!
			\!\!
			\!\!
			\!\!
			\max_{\theta}   \ \  &   \theta^\transpose \avgpsi      
			\label{e:ConvexApproxSima}
			\\
			\st  \ \  &  Q^\theta(z) \le  c (z)  +   \uQ^\theta ( \Dds(z) )  \,, \quad z\in \zstate
			\label{e:ConvexApproxSimb}
		\end{align} 
		This is a convex program of dimension $d$.
		
		\label{e:ConvexApproxSim}%
	\end{subequations}%
	
	\begin{subequations}%

		Many model-free algorithms might be used to approximate a solution to \eqref{e:ConvexApproxSim}.   This paper focuses on the simplest instance, in which an approximation of the inequality constraint in \eqref{e:ConvexApproxSimb} is defined by $\Gamma_\Nsam(\theta) \le 0$, where $N$ is the time horizon, and
		\begin{align}
			\Gamma_\Nsam(\theta) &\eqdef  \frac{1}{N} \sum_{k=0}^{N-1} [\Tdiff_{k+1}^\circ(\theta)]_-  \,,   \qquad 
			\\
			[\Tdiff_{k+1}^\circ(\theta)]_-  &\eqdef \max \{ 0,  -  \Tdiff_{k+1}^\circ(\theta) \} 
		\end{align}
		with $\Tdiff^\circ_{k+1}(\theta)$ defined in \eqref{e:TDonline}.	

\label{e:cvxQloss}		
\end{subequations}%

	In addition to the loss function \eqref{e:cvxQloss},   the algorithm introduced next requires a convex regularizer $ \regBCQ_\Nsam(Q, \theta)$,  	penalty parameter $\kappa\ge 0$,    tolerance $\Tol\ge 0$,  and
a probability measure $\upmu$ on $\clB(\zstate)$;    this will be chosen to be discrete, and in some cases based on observed input-state pairs.

			\begin{subequations}
\wham{Convex Q-learning}%
			Given the data $\{z(k) : k\le N\}$,   solve 
			\begin{align} 
				\theta^\ocp  \in \argmin_\theta
				\ \  &  \bigl \{- \theta^\transpose \avgpsi   + \kappa \regBCQ_\Nsam (Q^\theta, \theta)\bigr\} 
				\\
				\st  \ \   
				& \Gamma_\Nsam(\theta)  \le \Tol 
				\label{e:LPQ1constraint}
			\end{align}
			
Slater's condition holds provided $\Tol>0$.
			\label{e:LPQ1reg}%
		\end{subequations}%  

\notes{primal dual needed for journal}
%\bl{Consider the primal-dual algorithm to solve \eqref{e:LPQ1reg}:
%	\begin{align*}
%		\theta_{n+1} &= \argmin_\theta
%		\bigl \{ - \theta^\transpose \avgpsi  + \lambda_n [\Gamma_\Nsam(\theta) - \Tol ] +\kappa \regBCQ_\Nsam (Q^\theta, \theta)  \bigr \}
%		\\
%		\lambda_{n+1} &=  \max \bigl\{  0,   \lambda_n + \beta_n[\Gamma_\Nsam (\theta_{n+1}) - \Tol]   \bigr\}
%	\end{align*}
%	The introduction of $\Tol > 0$ here is crucial: $\Gamma_\Nsam(\theta_{n+1}) \ge 0$ for any $n$, so $\{ \lambda_n \}$ is a non-decreasing sequence if $\Tol = 0$.}
	
	\begin{figure*}
		\centering
		\includegraphics[width=1\hsize]{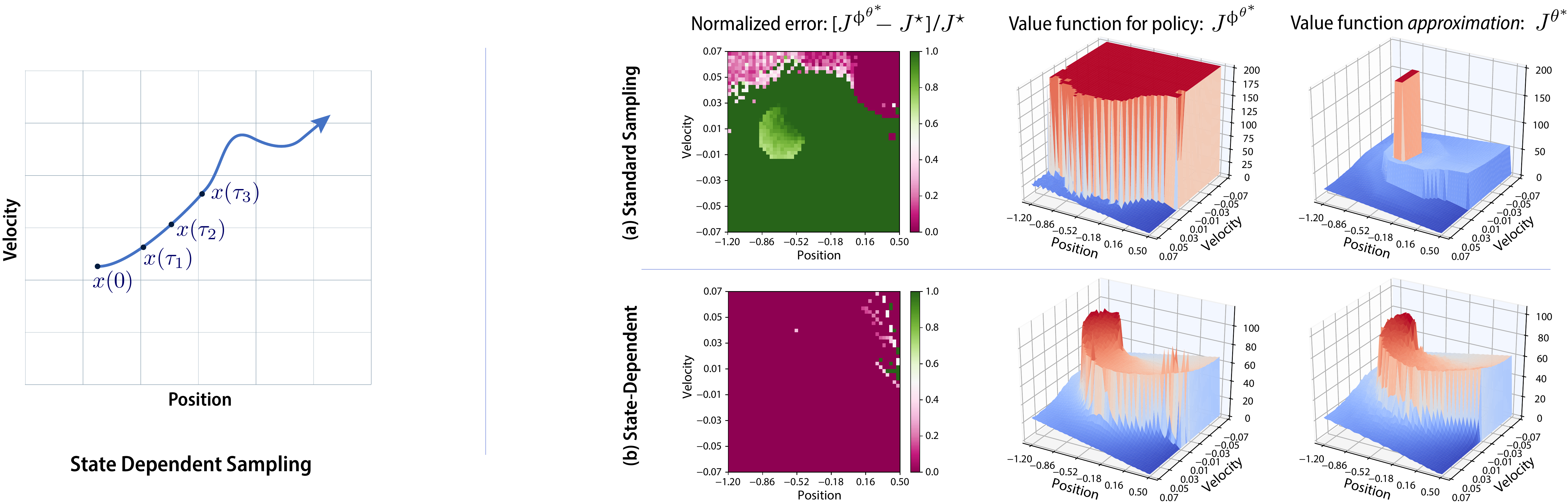}
		\caption{The figure on the left shows how the sampling times $\tau_k$ are chosen based on binning of the state space.  
			The right hand side shows value functions and their approximations for the mountain car example.    The first row illustrates failure of the algorithm due to numerical instability with fast sampling.   The second row shows that state-dependent sampling resolves this issue:  the value function approximation is very close to $Q^\oc$.  }
		\label[figure]{f:expNumInstab}
	\end{figure*}	
	\subsection{Exploration and  constraint geometry}
	\label{s:FanTheory}
	
	In this subsection only we impose an ergodicity assumption to ease analysis.   Using the compact notation, 
	$\csub{k} = c(z(k)) $   and $ \psisub{k} = \psi(z(k)) $ for $ k\ge 0 $,
%	\begin{equation}
%		\csub{k} = c(z(k)) \,,\qquad \psisub{k} = \psi(z(k)) \,,\qquad k\ge 0 \,,
%		\label{e:cpsisub}
%	\end{equation}
the following limits are assumed to exist,
	\begin{align*}
		\barpsi \eqdef \lim_{N\to\infty} \frac{1}{N} \sum_{k=0}^{N-1}\psisub{k}  
		\qquad
		R^\psi   \eqdef \lim_{N\to\infty}  \frac{1}{N} \sum_{k=0}^{N-1}\psisub{k} \psisub{k}^\transpose
	\end{align*} 
	and the covariance is denoted $	\Sigma^\psi\eqdef  R^\psi - \barpsi\barpsi^\transpose$.     These definitions appear in TD-learning;  in particular, it is common to say that  $\{\psi_i\}$ are  \textit{linearly independent}   if   $R^\psi$ is full rank.   The stronger condition 	$\Sigma^\psi>0$ is imposed in the following:  
	\begin{theorem}
		\label[theorem]{t:FanTheorem2}
		The   constraint region \eqref{e:LPQ1constraint} is always non-empty since it contains $\theta=\Zero$. 
		If $\Sigma^\psi$ is full rank, then the constraint region is bounded for all $N\ge 1$ sufficiently large.  
	\end{theorem}

The   lemmas that follow will quickly imply the conclusion of the theorem.   	  
	Proofs are contained in the Appendix.
	\begin{lemma}
		\label[lemma]{t:FanTheorem}
		Suppose that  the constraint region \eqref{e:LPQ1constraint} is unbounded for some $N\ge 1$. Then there exists a non-zero vector $\ctheta\in\Re^d$ such that $Q^{\ctheta}(z(k))$ is non-decreasing:
		\begin{equation}
			Q^{\ctheta}(z(k))  \ge Q^{\ctheta}(z(k-1))  \,,   \qquad    1\le k\le N
			\label{e:Qincreasing}
		\end{equation}
	%	\qed
	\end{lemma}

	\begin{lemma}
		\label[lemma]{t:increasingQimplication}
		Suppose that  \eqref{e:Qincreasing} holds for a fixed non-zero parameter $\ctheta\in\Re^d$,  and 
		every $N$ and every $1\le k\le N$. 
		Then
		$\Sigma^\psi$ is not full rank. %\qed
	\end{lemma}

	\wham{Proof of \Cref{t:FanTheorem2}} 
	We prove the contrapositive:  $\neg B \Longrightarrow \neg A$ where $A$ represents the logical expression ``\textit{$\Sigma^\psi$ is full rank}'',  and $B$ represents ``\textit{the constraint region  \eqref{e:LPQ1constraint}  is bounded for all $N\ge 1$ sufficiently large}''.

	The constraint region  \eqref{e:LPQ1constraint}  is non-decreasing with $N$, so that under  $\neg B$ it follows that the constraint region is unbounded for \textit{every $N$}.   \Cref{t:FanTheorem} then implies that the assumptions of  
	\Cref{t:increasingQimplication} hold:  this lemma tells us that  $\Sigma^\psi$ is not full rank, which is $\neg A$.
	\qed
	
	\subsection{Duality}
	\label{sec:dual}
	We consider the dual of \eqref{e:LPQ1reg} in the limiting case where $\Tol =0$, and $\kappa=0$, giving
	\begin{subequations}
		\begin{align*} 
			\theta^\ocp = \argmax_\theta
			\ \  &   \theta^\transpose \avgpsi \qquad
			\st  \   -\Tdiff^\circ_{k}(\theta)  \le 0 \,, \quad 1 \le k \le N
			\label{e:LPQ1constraintLim}
		\end{align*}%    
		\label{e:LPQ1regLim}%
	\end{subequations}% 
	The constraints are convex, but not linear.   An LP is obtained through the equivalent representation of the constraints:
	\begin{equation}
		\begin{aligned}
			Q^\theta(z(k-1)) - \csub{k-1} - Q^\theta(x(k), u) \le 0
		\end{aligned}
		\label{e:newconstt}
	\end{equation}
	for each $1\le k \le N$ and $u\in\ustate(x)$.
	
	For simplicity, in this section only we take $\ustate(x) =\ustate$ for each $x$.
	Denote $u^i$ the $i^{\text{th}}$ element in $\ustate$ for $1\le i \le \nU$, and $\avgpsi\in\Re^d$ is defined in \eqref{e:mupsi}.

	A column vector $\clC$ of dimension $n_C =\nU\times N$ and matrix $\clA$ of dimension $d\times  n_C$ are defined as follows:
	\begin{align*}
		\clC &\eqdef [C , \cdots , C]^\transpose \,,  
		&& 
		C \eqdef [\csub{0}, \csub{1}, \hdots, \csub{N-1}] 
		\\
		\clA &\eqdef
		\begin{bmatrix}
			A_1\\
			A_2
			\\
			\vdots\\
			A_{\nU}
		\end{bmatrix} \,, 
		&&  A_i = \begin{bmatrix}
			(\psisub{0} - \psi(x(0), u^i))^\transpose, \\
			(\psisub{1} - \psi(x(1), u^i))^\transpose, \\
			\vdots
			\\
			(\psisub{N-1} - \psi(x(N), u^i))^\transpose \\
		\end{bmatrix}
	\end{align*}
	where $1\le i\le \nU$.
	Then we arrive at the following LP:
	
	%\begin{subequations}
	\begin{align}
		\max_{\theta \in \Re^d} \quad &   \theta^\transpose \avgpsi
		\qquad
		\st \ \  \clA \theta \le \clC
		\label{e:LPPrimal}
	\end{align}
	%\label{e:LPPrimal}
	%\end{subequations}
	
	See  \cite[Ch.~4]{lue03}  for the derivation of its dual:
	\begin{subequations}
		\begin{align}
			\min \ \  
			& \sum_{k=1}^{N} \sum_{u\in \ustate} \varpi_{k, u} \csub{k-1}
			\label{e:LPDualNewObj}
			\\
			\st \ \  
			& \sum_{k=1}^{N} \sum_{u\in \ustate} \varpi_{k, u}\{  \psisub{k-1} - \psi(x(k), u) \} = \langle \upmu, \psi \rangle
			\label{e:newdualconst} 
		\end{align}
		where the minimum is over all  $\varpi \in \Re^{N\times \nU}$ satisfying $\varpi_{k, u} \ge 0 $ for each $ 1 \le k \le N, u\in\ustate$.
		\label{e:LPDualNew}
	\end{subequations}
	
	If $\varpi^\ocp$ is an optimizer of \eqref{e:LPDualNew}  and  $\theta^\ocp$   an optimizer of \eqref{e:LPPrimal}, then complementary slackness is expressed
	\begin{align}
		\varpi_{k, u}^\ocp [-Q^{\theta^*}(z(k-1)) + \csub{k-1}  + Q^{\theta^\ocp}(x(k), u)] = 0
		\label{e:comslack}
	\end{align}
	with $1 \le k \le N, u\in\ustate$.
	\Cref{t:FanPrimalDual} summarizes an immediate but interesting consequence.

	\begin{proposition}
		\label[proposition]{t:FanPrimalDual}
		Suppose that $(\theta^\ocp,\varpi^\ocp)$ are primal-dual optimizers.  
		If  $\varpi_{k^\circ, u^\circ}^\ocp>0$ for some $k^\circ$ and $u^\circ\in\ustate$ then the following holds:
		\begin{equation}
			\begin{aligned}
				\!\!
				\!\!
				\!\!
				\!\!
				0=
				\min_u &\{
				-Q^{\theta^\ocp}(z(k^\circ-1)) + \csub{k^\circ-1}  + Q^{\theta^\ocp}(x(k^\circ), u))  \}  
				\\
				&
				u^\circ\in \argmin_u   Q^{\theta^\ocp}(x(k^\circ), u))   
			\end{aligned} 
			\label{e:TCOE_PrimalDual} 
		\end{equation}
	\end{proposition}

	\wham{Proof}
	We have by feasibility of $\theta^\ocp$, for every $u\in\ustate(x)$,
	\begin{align*}
		-Q^{\theta^\ocp}(z(k^\circ-1)) + \csub{k^\circ-1}  + Q^{\theta^\ocp}(x(k^\circ), u) \ge 0
	\end{align*}
	and if $\varpi_{k^\circ, u^\circ}^\ocp>0$ then  \eqref{e:comslack} implies that we achieve this lower bound:
	\[
	-Q^{\theta^\ocp}(z(k^\circ-1)) + \csub{k^\circ-1}  + Q^{\theta^\ocp}(x(k^\circ), u^\circ)) = 0
	\]
	This establishes  the desired conclusion.
	\qed
	
	\begin{figure*}[h]
		\centering
		\includegraphics[width=1\hsize]{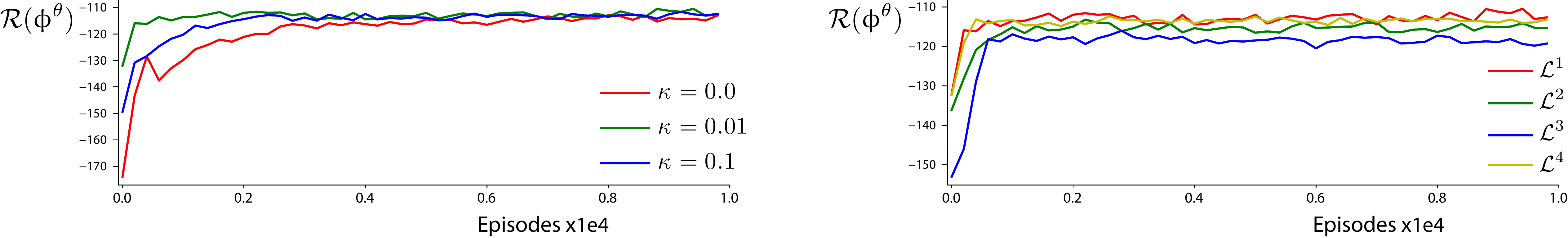}
		\caption{Performance comparison for Mountain Car with different $\kappa$ and different regularizers.}
		\label[figure]{f:RegComp}
	\end{figure*}	
	We have a clearer interpretation of the dual variable in the tabular setting, based on the following representation for the solution to the primal.
	Let $\{ z^\oc(k) = (x^\oc(k), u^\oc(k)) : k\ge 1\}$ be an optimal solution obtained with $z^\oc(0) $ chosen randomly according to $\upmu$.    Then,  
	\begin{equation}
		\langle \upmu, Q^\oc \rangle  =   \Expect\Bigl[  \sum_{k=0}^\infty c(x^\oc(k), u^\oc(k))   \Bigr]
		\label{e:mu-Q*}
	\end{equation}
	In the tabular setting, the basis is a collection of indicator functions:
	\[
	\psi_i(z) = \ind\{ z = z^i \}  \,,\qquad 1\le i\le d
	\]
	%FL2SM: "so that $d = n\times \nU -1$.": the lowercase $n$ comes out of nowhere. It should be changed to cardinality of \state.
	where $\{z^i : 1\le i\le d \} =  \state\times\ustate \setminus \{ z^e\}$, so that $d = |\state|\times \nU -1$.
	The equilibrium is omitted since we know that  $Q^\oc(z^e)=0$.
	
%	The proof of \Cref{t:PrimalDualOccMea} is postponed to the Appendix.

	\begin{proposition}
		\label[proposition]{t:PrimalDualOccMea}
		Consider the tabular setting,  and suppose that each state action pair in $ \state\times\ustate \setminus \{ z^e\}$ is visited at least once before time $N$.  Suppose also that the optimal policy $\fee^\ocp$ is unique.  Then   $Q^\oc$ is an optimizer of \eqref{e:LPPrimal},
		and  the dual variable   has the representation
		\begin{equation}
			\sum_{k=1}^N 
			\varpi^\ocp_{k,u} =    \sum_{k=0}^{N-1} \Prob\{  u^\oc(k ) = u\}   \,,   \quad u\in\ustate(x^\oc(k))
			\label{e:varpiTabular}
		\end{equation} 
%		\qed
	\end{proposition}

	\section{Numerical Results}
	\label{s:nres}
	
We survey here results from experiments on three examples from OpenAI Gym,
	Mountain Car,  CartPole, and Acrobot,
	focusing on four topics:
	(i) State-dependent sampling to improve algorithm performance;
	(ii) balancing the trade-off between exploration and exploitation;
	(iii) stability and consistency of convex Q-learning across different domains.
We also consider LQR to demonstrate that convex Q-learning is successful in cases where standard Q-learning diverges.
	
% No need for CDC:	Mountain Car \cite{moo90}, CartPole \cite{barsutand83a} and Acrobot \cite{geramifard2015rlpy}, 
	
	The Q-function approximations were defined by a linear function class in each experiment: $\{ Q^\theta = \theta^\transpose \psi: \theta\in\Re^d \}$,
	in which the basis functions 
	took the following separable form:
	\begin{align}
		\psi_{i,j}(z)= \begin{cases}
			0 \,, \quad  \text{if $z = z^e$},
			\\
			k_i(x) \ind\{ u = u^j \}    \,, \quad 1\le i \le d_x \,, \text{else}
		\end{cases}
		\label{e:stApp}
	\end{align}
	The functions $\{k_i\}$ were obtained using  the Python  function 
	{\small\tt sklearn.kernel\_approximation.RBFSampler} with $d_x = 250$;  see \cite{rahimi2007random, rahimi2008weighted}.
	Parameters for this function were chosen to be
	\[
	\sigma =  [0.05, 0.49, 0.93, 1.37, 1.81, 2.24, 2.68, 3.12, 3.56, 4.00]
	\]
	% FL2SM: d_x is not defined. Change made: I define d_x in \cref{e:stApp}.
	The dimension of $\theta$ is thus $d= d_x\times \nU$.

	Other approaches were investigated, such  as tile coding \cite{miller1990cmas}, radial basis functions, and binning.  We omit results using these approaches since the results were less reliable for the same dimension.	
	%\begin{wrapfigure}[14]{R}[0pt]{0.4\textwidth}
	%		%\vspace{-2em}
	%	\centering 
	%	\includegraphics[width=0.95\hsize]{figures/samTechDetails.pdf}
	%	\caption{State-dependent sampling for Mountain Car }
	%\end{wrapfigure}
	
	%\begin{figure}
	%\label[figure]{f:subsamp}
	%	%\vspace{-2em}
	%\centering 
	%\includegraphics[width=0.75\hsize]{figures/samTechDetails.pdf}
	%\caption{State-dependent sampling for Mountain Car }
	%\end{figure}

	\wham{Numerical Instability due to Fast sampling}  

	We frequently observe, especially when using a basis obtained through binning, that $\Tdiff_{k+1}^\circ (\theta) \approx \csub{k} \ge 0$, only because $\| x(k+1) - x(k) \|$ is small, for which the  constraint $\Tdiff_{k+1}^\circ(\theta) \ge 0$ is vacuous. This is purely an artifact of fast sampling: for example, the sampling interval $\Delta$ chosen in Mountain Car is $10^{-3}$. 
	Increasing the sampling interval will address this numerical challenge, but create new challenges because of the introduction of delay.
	
 	It is demonstrated here that state-dependent sampling can be designed to address this numerical challenge.
	
The sampling scheme begins with binning:  express the state space as a disjoint union  $\state = \cup_i \state_i$, and select  sampling times 
	$\{\tau_k\}$ so that sampled states are in distinct bins as illustrated in the plot on the left hand side of  \Cref{f:expNumInstab}.    
	The sampling times are defined recursively as follows:  choose an upper limit $\bar n$, take $\tau_0 = 0$, and for all $k \ge 0$ denote 
	\begin{equation}
		\begin{aligned}
			\tau_{k+1} &= \min\{ \tau_{k} + \bar n, \tau^0_{k+1} \}  
			\\
			\tau^0_{k+1} &= \min\{ j\ge \tau_{k} + 1: \text{Bin}(x(\tau_{j})) \not= \text{Bin}(x(\tau_{k})) \}
		\end{aligned}
		\label{e:tauk}
	\end{equation}
	where $\text{Bin}(x)$ denotes the index for the bin containing $x$.
	It is assumed  that the input takes a constant value on the interval $[\tau_{k}, \tau_{k+1})$ for each $k$, which justifies the the introduction of the cumulative cost, 
	%FL2SM: use macro \varC for state-dependent cost.
	\begin{align*}
		\varC_{\tau_{k}} = \sum_{j = \tau_{k}}^{\tau_{k+1} - 1}c(x(j), u(\tau_{k}))
	\end{align*}
	The temporal difference sequence is then redefined:
	\begin{align}
		\sTdiff_{k+1}^\circ(\theta) &\eqdef -Q^\theta(z(\tau_{k})) + \varC_{\tau_{k}} + \uQ^\theta(x(\tau_{k+1})) 
		\label{e:TdiffRedefined}
	\end{align}

	\begin{figure*}[h]
		\centering
		\includegraphics[width=\hsize]{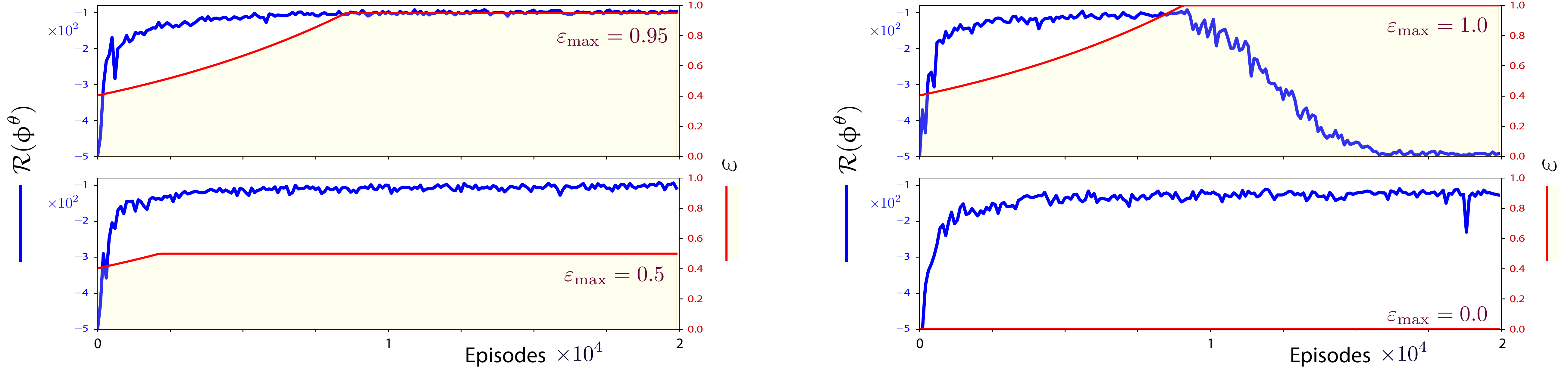} 
		\caption{Comparison of average cumulative rewards for the Acrobot with four different values of $\epsy_{\text{max}}$}
		\label[figure]{f:BadExp}
	\end{figure*}		
	\wham{Episodic Convex Q Learning}
	
	The basic algorithm \eqref{e:LPQ1reg} has been the focus of the previous section mainly for ease of analysis.   The experiments that follow are designed to approximate the solution to the finite time-horizon optimal control problem.
		
	\begin{subequations} 
		
In each example there is a goal set $\state^E \subset \state$ after which the state is reset. For training, we restart when the goal is reached. 
The initial condition $x^n(\Treset_{n})$ for  the $n$th episode after restart is chosen uniformly at random from $\state^\circ \subset \state$.    The successive restart times are defined by $T_0=0$, and   
		\begin{align*}
			\Treset_{n+1} &\eqdef \min\{  \bar n_E,   \Treset^\circ_{n+1} \}
			\\
			\Treset^\circ_n &\eqdef \min \{ \tau_k \ge T_n  : x^n(\tau_k) \in \state^E \}\,, \quad n\ge 0\,,
\end{align*} 
with $\bar n_E$ an upper limit imposed on the episode length, and $\{x^n(\tau_k) \}$ the trajectory from the $n$th episode.
 
		With  $B_{n+1} \eqdef T_{n+1} - T_n$, the parameter estimates are updated only at these times according to   
		\begin{align*} 
			\!\!\!\!\!
			\theta_{n+1} = \  & \argmin  \Bigl\{   -   \theta^\transpose \avgpsi       +  	\kappa\regBCQ_n ( \theta) 	
			+ \half \frac{1}{\alpha_{n+1}  }\| \theta - \theta_n \|^2      \Bigr\}  
			\\
			&\st  \ \    \Gamma_{\Treset_n}(\theta) =\frac{1}{B_{n+1}} \sum_{k=\Treset_{n}}^{\Treset_{n+1} - 1} [\sTdiff_{k+1}^\circ(\theta)]_-  \le \Tol
			%\label{e:theup}
		\end{align*} 
in which $\{ \alpha_{n+1} \}$ plays a role similar to a step-size sequence.  A constant value worked well in all experiments.

	\end{subequations}

%FL2SM: revision:  the regularizer $\regBCQ_n$ taken to be $	\regBCQ_n^1$   in \Cref{f:reg}.

\begin{figure}[h]
		\centering
		\includegraphics[width=0.9\hsize]{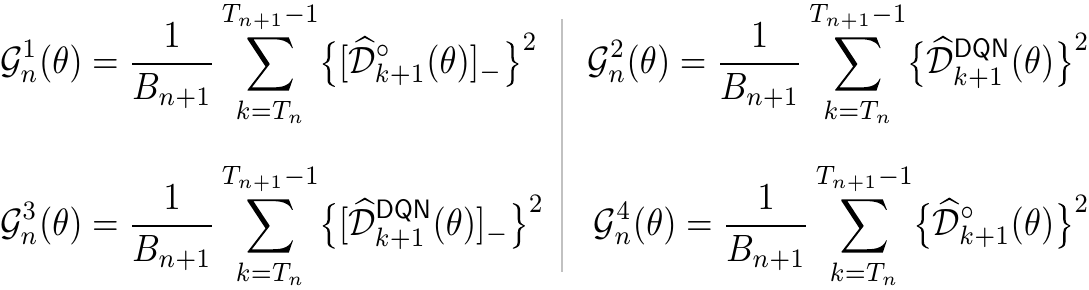} 
		\caption{Four regularizers considered in convex Q experiments.}
		\label[figure]{f:reg}
	\end{figure}

%\rd{revise:}	
The sensitivity of performance with respect to the coefficient $\kappa$ was investigated for each regularizer.    
The plots on the left hand side of \Cref{f:RegComp} were obtained using the regularizer $\regBCQ_n = 	\regBCQ_n^1$, defined in \Cref{f:reg}.    The best value in these experiments is  $\kappa = 0.01$.

	With  $\kappa = 0.01$ fixed, we then  compared performance of convex Q-learning using the regularizers  shown in \Cref{f:reg},  with $
	\sTdiff_{k+1}^\circ(\theta)$   defined in \eqref{e:TdiffRedefined},   and 		
	%FL2SM: Typo fixed: c(z(\tau_{k}  ) should be c_\tau{k}
		$\sTdiffDQN_{k+1} (\theta) \eqdef   -  Q^\theta (z(\tau_{k}) )   +   \varC_{\tau_{k}}  +   \uQ^{\theta_n} (x(\tau_{k+1})) $  (note dependency on $n$).

	%% Is this resolved?  Asking on September 14
%\notes{FL2SM: 
%Definition of $\Tdiff_{k+1}^{\text{DQN}} $ is a little bit odd here if we want to save examples of regularizer  can be found in Figure out a way so that \Cref{f:reg}. \Cref{f:reg} is be self-contained.
%	\\
%	\bl{Ignore my comment in draft sent to our co-authors!   I didn't look closely.   I think it is fine} 
%	}
	
	%\begin{equation*}
	%	\begin{aligned}
		%		\regBCQ^1_n (\theta) &= 
		%		\frac{1}{B_{n+1}} \sum_{k=\Treset_{n}}^{\Treset_{n+1} - 1} \bigl\{[\sTdiff_{k+1}^\circ(\theta)]_- \bigr\}^2 
		%		\\
		%		\regBCQ^2_n (\theta) &= 
		%		\frac{1}{B_{n+1}} \sum_{k=\Treset_{n}}^{\Treset_{n+1} - 1} \bigl\{[\sTdiffDQN_{k+1}(\theta)] \bigr\}^2 
		%		\\
		%		\regBCQ^3_n (\theta) &= 
		%		\frac{1}{B_{n+1}} \sum_{k=\Treset_{n}}^{\Treset_{n+1} - 1} \bigl\{[\sTdiffDQN_{k+1}(\theta)]_- \bigr\}^2	
		%		\\
		%		\regBCQ^4_n (\theta) &= 
		%		\frac{1}{B_{n+1}} \sum_{k=\Treset_{n}}^{\Treset_{n+1} - 1} \bigl\{[\sTdiff_{k+1}^\circ(\theta)] \bigr\}^2
		%	\end{aligned}
	%%\label{e:fourreg}
	%\end{equation*}	

The plots in \Cref{f:RegComp} show that $\regBCQ^1_n$ gives the best performance;
	this regularizer was chosen in  all subsequent experiments, with $\kappa = 0.01$.

	\wham{Exploration}
	
	The numerical results expose one gap in current theory:
	we obtain a convex program only when the input does not depend on parameters.      
	In each experiment, performance of convex Q-learning was greatly improved with an epsilon-greedy policy.    The following experiments are based on the following input sequence for training:
	%FL2SM: change $\theta_k$ to $\theta_n$ since \theta and \epsy are updated episodically.
	\begin{align*}
		\feex^{\theta_n}(u \mid  x)
		&\eqdef    \Prob\{u(\tau_k) = u \mid x(\tau_k) = x\} 
		\\
		&=     (1 - \epsy_n) P_E(u)   +     \epsy_n  \ind\{ u = \fee^{\theta_n}(x)  \}
	\end{align*}
	where $P_E$ is the uniform distribution on $\ustate(x)$.    
	The exploration parameter was chosen to be monotone increasing in $n$:  for parameters $ \xi ,  \epsy_{\text{max}} $  chosen in the interval $[0,1]$, 
	\begin{align}
		\epsy_{n+1} = \max \{ (1 +\xi)\epsy_n  , \epsy_{\text{max}}  \}  \,,\quad n\ge 0
		\label{eq:decayeps}
	\end{align}
	initialized with $\epsy_0>0$ (a small constant).
	
	\wham{Validation}   A parameter $\theta\in\Re^d$ was evaluated by conducting  $N$ independent experiments under the $Q^\theta$-greedy policy.    The $n$th experiment 
	%FL2SM: change $\{ x^n(k), u^n(k) ... to $\{ z^n(\tau_k)
	results in a time $\Treset_n$ at which the terminal state is reached,  and data $\{ z^n(\tau_k):  \Treset_{n}\le k\le \Treset_{n+1} \}$. 
	
	The examples from OpenAI Gym are based on a  one-step reward function $r$;  $c=-r$ was used in the algorithms described above, but for validation we computed the average cumulative reward:
	\begin{align}
		\mathcal{R}(\fee^\theta) \eqdef  -\frac{1}{N} \sum_{n=1}^{N} \bigg [\sum_{k=\Treset_{n}}^{\Treset_{n+1} - 1} \varC_{\tau_{k}} \bigg ]  
		\label{e:cummReward}
	\end{align}
	%FL2SM: within the interval [\tau_k, \tau_{k+1}), u^n(\tau_k) takes the same value.
	with $u^n(\tau_k) = \fee^\theta(x^n(\tau_k)) $.

	\wham{Experimental results in three examples}   To test consistency of outcomes we performed repeated runs in several control system examples.  
	
	In each case,    $50$ independent experiments where executed,  in which  the initial condition $\theta_0$ was chosen independently according to a normal distribution $N(0, I)$, and initial conditions were chosen uniformly   from $\state^\circ$.
	
	Selected results are collected in \Cref{f:allresults}, which show that convex Q-learning   succeeded in solving the three examples, with remarkable consistency across runs.  	
	
	\Cref{f:expNumInstab} shows results obtained for the Mountain Car problem:  the first row shows the failure of the algorithm with fast sampling (time-steps from the standard model);  the value function approximation is very poor, and it was found that the resulting $Q^\theta$-greedy policy is unacceptable.  
	The second row shows nearly perfect approximation of the true value function $Q^\oc$ when using state-dependent sampling. 
	
	It will not surprise many readers to learn that the parameter $\epsy_{\text{max}} $ defined in \eqref{eq:decayeps} 
	should be chosen with some care.     \Cref{f:BadExp} shows results from four implementations of convex Q-learning for the  Acrobot:  each figure includes two plots as a function of episode $n$:
	the exploration parameter $\epsy_n$ and the cumulative reward $\clR(\fee^{\theta_n})$ obtained from parameter $\theta_n$  (definition in \eqref{e:cummReward}).  
	The plots shown on the left hand side use $\epsy_{\text{max}} =0.95$ and $\epsy_{\text{max}} =0.5$.  It is seen that the reward quickly reaches the desired value that is considered a solution to the Acrobot example: $\clR^\ocp \ge -100$.  	The plot on the lower right shows that there is greater bias with pure exploration, i.e., $\epsy_{\text{max}} =0$.

	The algorithm with  $\epsy_{\text{max}} =1$ fails:
The plots on the top right hand side show that performance drastically drops at episode 8,000,  when $\epsy_n$ reaches about $0.975$.    	
	
\wham{Comparison with standard Q-learning} 
A very simple example illustrates the striking difference between convex Q-learning and the standard algorithm  \eqref{e:Q}.

Consider the one-dimensional LQR model with dynamics $\dot x = u$, and quadratic cost $c(x, u ) = x^2 + u^2$; the optimal policy is linear state feedback,  and the Q-function is quadratic.   This motivates the basis 
$Q^\theta(x, u) \eqdef \theta^\transpose \psi(x, u)$ with $\psi(x, u) = [x^2, 2xu, u^2]^\transpose$. 

Convex Q-learning and the standard Q-learning algorithm were compared with the input for training:   $u(t) = 
  -K^* x(t) + \sum_{i=1}^{10} \sin((10+40v_i)t)$,   with $K^*$ the optimal gain and $v_i $ uniformly sampled from $ [-1, 1]$.   With $\Tol = 0.01$ and no regularizer, the solutions of convex Q-learning \eqref{e:ConvexApproxSim} are consistent after a short run.   
  
In contrast, estimates from the algorithm  \eqref{e:Q} are typically divergent---\Cref{f:compCvxQDQN} shows the evolution   of the second component of $\{\theta_n :  n\ge 0\}$ from two initial condtions.
  The source of instability for large initial condition $\theta_0$ is lack of Lipschitz continuity of $f$ and $\barf$:   the minimizer over $u$ of $Q^\theta(x,u)$ defines the $Q^\theta$-greedy policy, which in this case takes the explicit form $    \fee^\theta(x) = - K^\theta x$,   with $K^\theta = \theta(2)/\theta(3)$.  Consequently,   
  \[
  	\uQ^\theta(x) = Q^\theta(x, \fee^\theta(x) ) = \beta_\theta  x^2\,,\quad \beta_\theta\eqdef \theta(1) - \theta(2)^2/\theta(3)
  \]
and from  \eqref{e:Q} and \eqref{e:WatkinsRelax}

\[
\begin{aligned}
\barf(\theta)  &=  \Expect
			\bigl[      \Tdiff^\circ_{k+1}(\theta) \psi (z(k))    \bigr]  
   \\
       \Tdiff^\circ_{k+1}(\theta)      & =   - \theta^\transpose \psi(z(k))  + x(k)^2 + u(k)^2   + \beta_\theta x(k+1)^2
\end{aligned} 
\]
The mean flow $\ddt \odestate = \barf(\odestate)$ is not globally asymptotically stable [it has finite escape time for initial condition satisfying $\odestate_3(0)<0$ and large $\odestate_2(0)>0$].    And, even if the ODE were stable,  the fact that $\barf$ is not Lipschitz continuous violates an essential assumption in stochastic approximation.
  
  \notes{I may have expanded too much!
  \\
  Note that I corrected a bug: you had written $ \phi^\oc$ instead of $\fee^\theta$---both "oc" and "phi" aren't right!
  }

  %gives $K = 0.618$ and $M = 1.618$. 
%. The optimal stated-feed back gain $K^* = 0.618$ and the optimal value function is  $J^*(x) = x^\transpose M^* x$ with $M^* = 1.618$. Consider Both of the two algorithms are trained using data generated by the same randomized policy:
%	$

%	\clearpage

	\section{Conclusion and Future Work}  
	
	%\rd{Work in progress!}
	%	(i) The dual of convex Q-learning is not precisely  Manne's LP or a version of logistic Q-learning, but has similar structure that reveals the need for regularization to  avoid over-fitting.   
	%	(ii) Necessary and sufficient conditions for a bounded solution to the Q-learning LP with linear function approximation are obtained for the first time.   
	%	(iii)   Simulation studies reveal numerical challenges when addressing sampled-data systems based on a continuous time model.  The challenge is addressed using state-dependent sampling.   The theory is illustrated with applications to examples from OpenAI gym. 
	%	In this paper, we presented an approach, called convex Q-learning, that builds on the LP formulation of optimal control of Manne. 
	%	We studied the dual of this algorithm and revealed the need for regularization. We provided theoretical analysis on the conditions of a bounded solution to the Q-learning LP with linear function approximation, and demonstrated its performance on several experiments. 
	
	Convex Q-learning is a recent approach to reinforcement learning, whose main appeal is that we have a better understanding of what the algorithm is attempting to solve.    There are of course many open questions that will be explored in future research:   can we obtain performance bounds in non-ideal settings?   There is some hope, given that Lyapunov functions might be introduced in an augmented LP.    Can we extend convergence theory to the parameter-dependent policies used in the numerical results?  Can we obtain more efficient algorithms using deeper theory of quasi-stochastic approximation \cite{laumey22d,laumey22e}?
	 We are also considering alternate algorithm architectures for application of RKHS techniques.

	%\clearpage
	
	\normalem
	\bibliographystyle{abbrv}
	\bibliography{IEEEabrv,strings,markov,q,DPLP,kernel,FLextras}  
	%fansref
	
%	
%	{}
%\null  
%\null  %Needed with \usepackage{flushend}
%\end{document} 

%	\choreS{All done?  
%	\\
%	SM2All
%		\\
%		To do for final submission in the summer of 2022:
%		\\
%		1.  Use PR averaging in the numerics.   
%		\\
%		2.   See if we can prove a converse, that bounded feasibility set  implies  full rank covariance }
	
	\appendix
	
	\section{Appendix}

	\subsection{Theory behind   \Cref{t:FanTheorem2}}
	
	\wham{Proof of \Cref{t:FanTheorem}}
	If the constraint region $\eqref{e:LPQ1constraint}$ is unbounded for some $N \ge 1$, then for each $m\ge 1$, there exists $\theta^{m}$ such that $\| \theta^{m} \| \ge m$, and
	%\begin{align}
	%\label{e:intolp}
	%	\frac{1}{N} \sum_{k=1}^{N} \max\bigg\{0, Q^{\theta^{m}}(z(k-1)) - \csub{k-1} - \uQ^{\theta^{m}} (x(k))\bigg\} \le \Tol
	%\end{align}
	\begin{align}
		\label{e:intolp}
		\frac{1}{N} \sum_{k=1}^{N} \max\bigg\{0, \Tdiff^\circ_{k}(\theta^{m})\bigg\} \le \Tol
	\end{align}
	Dividing \eqref{e:intolp} by $\|\theta^{m}\|$ gives
	%\begin{equation}
	%\label{e:infineq}
	%	\frac{1}{N} \sum_{k=1}^{N} \max\bigg\{0, \frac{Q^{\theta^{m}}(z(k-1))}{\|\theta^{m}\|} - \frac{\csub{k-1}}{\|\theta^{m}\|} - \frac{\uQ^{\theta^{m}}(x(k) )}{\|\theta^{m}\|}\bigg\} \le \frac{\Tol}{\|\theta^{m}\|}
	%\end{equation}
	\begin{equation}
		\label{e:infineq}
		\frac{1}{N} \sum_{k=1}^{N} \max\bigg\{0, \frac{\Tdiff^\circ_{k}(\theta^{m})}{\|\theta^{m}\|}\bigg\} \le \frac{\Tol}{\|\theta^{m}\|}
	\end{equation}
	Denote $\ctheta^m =  {\theta^m}/{\| \theta^m\|}$. By the definition of $\uQ^{\theta^m}$, we have
	\begin{align*}
		\frac{1}{\| \theta^m \|} \uQ^{\theta^m}(x(k)) &=  \min_u  \Bigl\{ \frac{1}{\| \theta^m \|} Q^{\theta^m}(x(k), u)  \Bigr\}
		\\
		&= \uQ^{\ctheta^m}(x(k))
	\end{align*}
	Thus, we can write \eqref{e:infineq}  as:
	\begin{equation}
		\label{e:infeqde}
		\frac{1}{N} \sum_{k=1}^{N} \max\bigg\{0, \widehat{\Tdiff}_k^\circ(\ctheta^{m})\bigg\} \le \frac{\Tol}{\|\theta^{m}\|}
	\end{equation}
	with $\widehat{\Tdiff}_k^\circ(\theta^{m}) = Q^{\theta^{m}}(z(k-1)) - \frac{\csub{k-1}}{\|\theta^{m}\|} - \uQ^{\theta^{m}}(x(k) )$.
	Since $\|\ctheta^m\|  = 1$ for each $m$, there exists a convergent subsequence $\{ \theta^{m_i} \}$ with limit satisfying  $\| \ctheta  \| = 1$:
	\[
	\lim_{i\to\infty} \frac{ \theta^{m_i} }{ \| \theta^{m_i} \| } = \lim_{i\to\infty} \ctheta^{m_i} = \ctheta
	\]
	The inequality  \eqref{e:infeqde} then gives
	\begin{equation}
		\label{e:inq3}
		\begin{aligned}
			\frac{1}{N} \sum_{k=1}^{N} & \max\bigg\{0, Q^{\ctheta}(z(k-1)) - \uQ^{\ctheta}(x(k) )\bigg\} 
			\\
			&=\lim_{i\to\infty}\frac{1}{N} \sum_{k=1}^{N} \max\bigg\{0, \Tdiff(\ctheta^{m_i})\bigg\}
			\le 0
		\end{aligned}
	\end{equation}
	where $\widehat{\Tdiff}_k^\circ(\ctheta^{m_i}) = Q^{\ctheta^{m_i}}(z(k-1)) - \frac{\csub{k-1}}{\|\theta^{m_i}\|} - \uQ^{\ctheta^{m_i}}(x(k) )$ for each $i$.
	Each term in the sum on the left hand side of \eqref{e:inq3} must be zero, and hence  
	\[
	Q^{\ctheta}(z(k-1)) \le \uQ^{\ctheta}(x(k)) = \min_u Q^{\ctheta}(x(k), u)\le Q^{\ctheta}(z(k))
	\]
	for $1\le k \le N$, which is the desired conclusion.
	\qed

	\wham{Proof of \Cref{t:increasingQimplication}}
	%\wham{Proof of \Cref{t:increasingQimplication}}

	We write the assumed sequence of inequalities in the equivalent form
	\[
	Q^{\ctheta}(z(k))  - Q^{\ctheta}(z(k-1))  = \epsy_k \,,   \qquad   \textit{for all $k\ge 1$}
	\]
	where $\{\epsy_k : k\ge 1\}$ is a non-negative sequence.   We then have for any $N\ge 1$,
	\[
	\sum_{k=1}^{N} \epsy_k  = Q^{\ctheta}(z(N))  - Q^{\ctheta}(z(0)) 
	\]
	This implies that the sequence $\{\epsy_k\}$ is summable under the assumption that $\{z(k)\}$ is bounded.
	For two integers $1\le N_0< N$  write  
	\begin{align*}
		S_{N_0}^N \eqdef \sum_{k=N_0}^{N}\epsy_k  
		&=   Q^{\ctheta}(z(N))  - Q^{\ctheta}(z(N_0-1))   
		\\
		&= \{ \psisub{N} - \psisub{N_0-1} \}^\transpose \ctheta
	\end{align*}
	Multiplying each side on the left by  $ \psisub{N}$  gives the suggestive identity
	\[
	S_{N_0}^N  \psisub{N}
	=  \psisub{N}  \psisub{N}^\transpose \ctheta     -   \psisub{N} \psisub{N_0-1}^\transpose \ctheta
	\]
	Consequently,  for any fixed $N_0$,    averaging over $N$ gives,
	\begin{align*}
		\epsy_{N_0}^\psi  &\eqdef  \lim_{M \to\infty} \frac{1}{M }   \sum_{N=N_0}^{M}  S_{N_0}^N  \psisub{N} 
		\\
		&= \lim_{M \to\infty} \frac{1}{M }   \sum_{N=N_0}^{M} \psisub{N}  \psisub{N}^\transpose \ctheta
		\\
		&\qquad -  \lim_{M \to\infty} \frac{1}{M }\sum_{N=N_0}^{M} \psisub{N} \psisub{N_0-1}^\transpose \ctheta
		\\
		&= 
		R^\psi \ctheta     -   \barpsi  \psisub{N_0-1}^\transpose \ctheta
	\end{align*}
	The left hand side vanishes as $N_0\to\infty$. Averaging over $N_0$ and taking a limit completes the proof:
	\begin{equation}
		0 = 
		\lim_{M \to\infty} \frac{1}{M } 
		\sum_{N_0=1}^{M-1} 
		\epsy_{N_0}^\psi 
		=
		[R^\psi       -   \barpsi    \barpsi^\transpose ] \ctheta
		=\Sigma^\psi \ctheta  
		\label{e:NInf}
	\end{equation}
	\qed

	\subsection{Proof of \Cref{t:PrimalDualOccMea}}
	%\wham{Proof of \Cref{t:PrimalDualOccMea}}
	Fix any $ u^\circ\in\ustate$  and perturb the   constraints of \Cref{e:LPPrimal}  as follows:  
	\begin{equation}
		\begin{aligned}
			\eta(\epsy)   =    \max  \quad &\langle \upmu, Q \rangle
			\\
			\st   \quad &  Q(z(k-1)) - c_{k-1}- Q(x(k), u) 
			\\ 
			&\qquad \le \epsy \ind\{    u= u^\circ \}   \,,\quad 1\le k\le N
		\end{aligned}%
		\label{e:primmmmCorrected}%
	\end{equation}
	where the maximum is over all functions $Q\colon\state\times\ustate\to\Re$ satisfying $Q(x^e,u^e)=0$.
	Under the assumption of the proposition, it follows from Prop.~2.1 of \cite{mehmeyneulu21} that $\langle \upmu, Q^* \rangle \ge \langle \upmu, Q^\theta \rangle $ for any $\theta$. As a consequence, $Q^*$ is an optimizer of the LP \Cref{e:LPPrimal}, in the sense that  $Q^* = Q^{\theta^*}$ with $\theta^*_i = Q^*(z^i)$ for each $i$.
	The LP \eqref{e:primmmmCorrected} is precisely convex Q-learning with perturbed cost function:  $c_\epsy(z) = c(z) +  \epsy  \ind\{    u= u^\circ \} $ for $z=(x,u)\in\zstate$.
	
	Letting $\clC_\epsy$ denote the new constraint vector in the LP we have for all $\epsy$,
	\[
	\begin{aligned}
		\eta(\epsy)  &  =   \min_\varpi \max_\theta  \bigl(   \theta^\transpose  \avgpsi  +  \{ - \clA \theta + \clC_\epsy \}^\transpose \varpi \bigr)
		\\
		& \le \max_\theta  \bigl(   \theta^\transpose  \avgpsi +  \{ - \clA \theta + \clC_\epsy \}^\transpose \varpi^*  \bigr)  =   \clC_\epsy ^\transpose \varpi^*
	\end{aligned} 
	\]
	This has the sub-derivative interpretation,
	\begin{equation*}
		\begin{aligned}
			\eta(\epsy) \le \eta(0) + \{  \clC_\epsy -\clC \}^\transpose \varpi^*  =   \eta(0) + \epsy \sum_{k=1}^\infty \varpi^*_{k, u^\circ}
		\end{aligned}
	\end{equation*}
	To complete the proof we obtain an alternative representation for  $\eta(\epsy) $ that is affine in a neighborhood of $\epsy=0$.

	Uniqueness of $\fee^*$ implies that it remains optimal with this new cost function, for all   $|\epsy|$ is sufficiently small, so that
	\begin{equation}
		\begin{aligned}
			\eta(\epsy) &=    \Expect\Bigl[  \sum_{k=0}^\infty c_\epsy(x^*(k), u^*(k))   \Bigr] 
			\\
			&=  % \Expect\Bigl[  \sum_{k=0}^\infty c(x^*(k), u^*(k))   \Bigr]  
			\eta(0)   + \epsy \Expect\Bigl[  \sum_{k=0}^\infty\ind\{ u^*(k) = u^\circ \}   \Bigr]     
		\end{aligned}
		\label{e:aff}
	\end{equation}
	Hence $\eta$ is differentiable,  and the sub-derivative becomes a derivative, with
	\[
	\frac{d}{d\epsy} \eta(\epsy) \big|_{\epsy =0} =  \sum_{k=1}^\infty \varpi^*_{k, u^\circ}
	\]
	This combined with \Cref{e:aff} completes the proof.
	\qed

	{}
\null  
\null  %Needed with \usepackage{flushend}
\end{document}

%% file: bookmacrosRL_IEEE.tex
     \makeatletter
\newcommand\gobblepars{%
    \@ifnextchar\par%
        {\expandafter\gobblepars\@gobble}%
        {}}
\makeatother

\def\mindex#1{\index{#1}}

 \def\oc{\star}   % for optimal control
\def\ocp{*}   % for optimal parameter or function approximation

\def\Tdiff{\mathcal{D}}

\def\odestate{\upvartheta}

\def\Nsam{N}

\def\fee{\upphi}

\def\feex{{\breve{\fee}}}

\def\elig{\zeta}

\def\uQ{\underline{Q}}

\newcommand{\bbblot}{\raise1pt\hbox{\vrule height .4ex width .4ex depth .05ex}}
\long\def\defbox#1{\framebox[.9\hsize][c]{\parbox{.85\hsize}{%
\parindent=0pt
\baselineskip=12pt plus .1pt      % STYLE 
\parskip=6pt plus 1.5pt minus 1pt % CHANGES
 #1}}}

\long\def\beginbox#1\endbox{\subsection*{}%
\hbox{\hspace{.05\hsize}\defbox{\medskip#1\bigskip}}%
\subsection*{}}

\def\endbox{}

 \def\archival#1{} %Notes I'd like to save

\def\FRAC#1#2#3{\genfrac{}{}{}{#1}{#2}{#3}}

\def\ddt{{\mathchoice{\FRAC{1}{d}{dt}}%
{\FRAC{1}{d}{dt}}%
{\FRAC{3}{d}{dt}}%
{\FRAC{3}{d}{dt}}}}

\def\ddtp{{\mathchoice{\FRAC{1}{d^{\hbox to 2pt{\rm\tiny +\hss}}}{dt}}%
{\FRAC{1}{d^{\hbox to 2pt{\rm\tiny +\hss}}}{dt}}%
{\FRAC{3}{d^{\hbox to 2pt{\rm\tiny +\hss}}}{dt}}%
{\FRAC{3}{d^{\hbox to 2pt{\rm\tiny +\hss}}}{dt}}}}

\def\ddyp{{\mathchoice{\FRAC{1}{d^{\hbox to 2pt{\rm\tiny +\hss}}}{dy}}%
{\FRAC{1}{d^{\hbox to 2pt{\rm\tiny +\hss}}}{dy}}%
{\FRAC{3}{d^{\hbox to 2pt{\rm\tiny +\hss}}}{dy}}%
{\FRAC{3}{d^{\hbox to 2pt{\rm\tiny +\hss}}}{dy}}}}

\def\half{{\mathchoice{\FRAC{1}{1}{2}}%
{\FRAC{1}{1}{2}}%
{\FRAC{3}{1}{2}}%
{\FRAC{3}{1}{2}}}}

\newsavebox{\junk}
\savebox{\junk}[1.6mm]{\hbox{$|\!|\!|$}}

\def\argmin{\mathop{\rm arg{\,}min}}
\def\argmax{\mathop{\rm arg{\,}max}}

\def\Dds{\text{\rm F}}

\def\state{{\sf X}}

\def\ustate{{\sf U}} 
 
\def\zstate{{\sf Z}}

\def\bfmath#1{{\mathchoice{\mbox{\boldmath$#1$}}%
{\mbox{\boldmath$#1$}}%
{\mbox{\boldmath$\scriptstyle#1$}}%
{\mbox{\boldmath$\scriptscriptstyle#1$}}}}

\def\bfmY{\bfmath{Y}}

\def\bfmhhaY{\bfmath{\hhaY}} %\widehat{\widehat{Y}}}}
\def\bfmhhaY{\hbox to 0pt{$\widehat{\bfmY}$\hss}\widehat{\phantom{\raise 1.25pt\hbox{$\bfmY$}}}}

\def\clB{{\cal B}}
\def\clC{{\cal C}}
\def\clR{{\cal R}}

\def\clC{{\cal C}}

\def\eqdef{\mathbin{:=}}

\def\Prob{{\sf P}}

\def\Expect{{\sf E}}

\def\Zero{{\mathchoice{\mbox{\sf 0}}%  \lgmath
{\mbox{\sf 0}}%
{\mbox{\tiny \sf 0}}%
{\mbox{\tiny \sf 0}}}}

 \def\Tol{\text{\rm Tol}}

\def\ind{\bbbone}
  
 \def\epsy{\varepsilon}

\def\formtmp#1#2{{\vskip12pt\noindent\fboxsep=0pt\colorbox{#1}{\vbox{\vskip3pt\hbox to \textwidth{\hskip3pt\vbox{\raggedright\noindent\textbf{#2\vphantom{Qy}}}\hfill}\vspace*{3pt}}}\par\vskip2pt%
\noindent\kern0pt}}

\newenvironment{programcode}[1]{\ignorespaces\def\stmtopen##1{##1}%
\pagebreak[3]%
\formtmp{programcode}{#1}%
\endlinechar=-1\relax%
\nopagebreak[4]}{%
\noindent\textcolor{programcode}{\rule{\columnwidth}{1pt}}\vskip1pt\par\addvspace{\baselineskip}%
\endlinechar=13}
 
\def\barf{{\overline {f}}}

\def\barpsi{{\bar{\psi}}}

{\end{list}}

\def\ass(#1:#2){(#1\ref{#1:#2})}

\def\ritem#1{
\item[{\sf \ass(\current_model:#1)}]
}

\newenvironment{recall-ass}[1]{% 
\begin{description}
\def\current_model{#1}}{
\end{description}
}

\def\sq{\hbox{\rlap{$\sqcap$}$\sqcup$}}
\def\qed{\ifmmode\sq\else{\unskip\nobreak\hfil
\penalty50\hskip1em\null\nobreak\hfil\sq
\parfillskip=0pt\finalhyphendemerits=0\endgraf}\fi}

\newcommand{\blot}{\vrule height 1.1ex width .9ex depth -.1ex }
\def\qedb{\ifmmode\blot\else{\vspace{-.2cm}\unskip\nobreak\hfil
\penalty50\hskip1em\null\nobreak\hfil\blot
\parfillskip=0pt\finalhyphendemerits=0\endgraf}\fi}

\newcounter{rmnum}
\newenvironment{romannum}{\begin{list}{{\upshape (\roman{rmnum})}}{\usecounter{rmnum}
\setlength{\leftmargin}{8pt}
\setlength{\rightmargin}{6pt}
\setlength{\itemindent}{0pt}
}}{\end{list}}

\newcounter{anum}

\newcommand{\field}[1]{\mathbb{#1}}

\def\Re{\field{R}}

\def\Prob{{\sf P}}

\def\Expect{{\sf E}}

\def\transpose{{\intercal}}

\def\st{\text{\rm s.t.\,}}

\def\argmin{\mathop{\rm arg\, min}}
\def\ind{\hbox{\large \bf 1}}

\def\epsy{\varepsilon}

\def\haY{\widehat{Y}}

\def\hhaY{\hbox to 0pt{$\haY$\hss}\widehat{\phantom{\raise 1.25pt\hbox{Y}}}}

\def\haY{\widehat Y}

%% file: convexQLearning_09_14_EXTENDED_spm.bbl
\def\cprime{$'$}\def\cprime{$'$}
\begin{thebibliography}{10}

\bibitem{bai95}
L.~Baird.
\newblock Residual algorithms: Reinforcement learning with function
  approximation.
\newblock In A.~Prieditis and S.~Russell, editors, {\em Proc. Machine
  Learning}, pages 30--37. Morgan Kaufmann, San Francisco (CA), 1995.

\bibitem{bascurkraneu21}
J.~Bas~Serrano, S.~Curi, A.~Krause, and G.~Neu.
\newblock Logistic {Q}-learning.
\newblock In A.~Banerjee and K.~Fukumizu, editors, {\em Proc. of The Intl.
  Conference on Artificial Intelligence and Statistics}, volume 130, pages
  3610--3618, 13--15 Apr 2021.

\bibitem{bertsi96a}
D.~Bertsekas and J.~N. Tsitsiklis.
\newblock {\em Neuro-Dynamic Programming}.
\newblock Atena Scientific, Cambridge, Mass, 1996.

\bibitem{bor02a}
V.~S. Borkar.
\newblock Convex analytic methods in {M}arkov decision processes.
\newblock In {\em Handbook of {Markov} decision processes}, volume~40 of {\em
  Internat. Ser. Oper. Res. Management Sci.}, pages 347--375. Kluwer Acad.
  Publ., Boston, MA, 2002.

\bibitem{farroy06}
D.~P. de~Farias and B.~Van~Roy.
\newblock A cost-shaping linear program for average-cost approximate dynamic
  programming with performance guarantees.
\newblock {\em Math. Oper. Res.}, 31(3):597--620, 2006.

\bibitem{gor00}
G.~J. Gordon.
\newblock Reinforcement learning with function approximation converges to a
  region.
\newblock In {\em Proc. of the 13th Intl. Conference on Neural Information
  Processing Systems}, pages 996--1002, Cambridge, MA, 2000.

\bibitem{laumey22e}
C.~K. Lauand and S.~Meyn.
\newblock Approaching quartic convergence rates for quasi-stochastic
  approximation with application to gradient-free optimization.
\newblock {\em Proc. Conference on Neural Information Processing Systems
  ({NeurIPS})}, 2022.

\bibitem{laumey22d}
C.~K. Lauand and S.~Meyn.
\newblock Quasi-stochastic approximation: Design principles with applications
  to extremum seeking control.
\newblock {\em {IEEE} Control Systems Magazine (to appear)}, 2022.

\bibitem{leehe19b}
D.~Lee and N.~He.
\newblock Stochastic primal-dual {Q}-learning algorithm for discounted {MDPs}.
\newblock In {\em Proc. of the American Control Conf.}, pages 4897--4902, July
  2019.

\bibitem{leehe19}
D.~Lee and N.~He.
\newblock A unified switching system perspective and {ODE} analysis of
  {Q}-learning algorithms.
\newblock {\em arXiv}, page arXiv:1912.02270, 2019.

\bibitem{mehmeyneulu21}
F.~Lu, P.~G. Mehta, S.~P. Meyn, and G.~Neu.
\newblock Convex {Q}-learning.
\newblock In {\em American Control Conf.}, pages 4749--4756. IEEE, 2021.

\bibitem{lumehmeyneu22}
F.~Lu, P.~G. Mehta, S.~P. Meyn, and G.~Neu.
\newblock Convex analytic theory for convex {Q}-learning.
\newblock In {\em Conference on Decision and Control--to appear}, page~PP.
  IEEE, 2022.

\bibitem{lue03}
D.~Luenberger.
\newblock {\em Linear and nonlinear programming}.
\newblock Kluwer Academic Publishers, Norwell, MA, second edition, 2003.

\bibitem{maeszebhasut10}
H.~R. Maei, C.~Szepesv\'{a}ri, S.~Bhatnagar, and R.~S. Sutton.
\newblock Toward off-policy learning control with function approximation.
\newblock In {\em Proc. ICML}, pages 719--726, USA, 2010. Omnipress.

\bibitem{man60a}
A.~S. Manne.
\newblock Linear programming and sequential decisions.
\newblock {\em Management Sci.}, 6(3):259--267, 1960.

\bibitem{mehmey09a}
P.~G. Mehta and S.~P. Meyn.
\newblock Q-learning and {Pontryagin's} minimum principle.
\newblock In {\em Proc. of the Conf. on Dec. and Control}, pages 3598--3605,
  Dec. 2009.

\bibitem{melmeyrib08}
F.~S. Melo, S.~P. Meyn, and M.~I. Ribeiro.
\newblock An analysis of reinforcement learning with function approximation.
\newblock In {\em Proc. ICML}, pages 664--671, New York, NY, 2008.

\bibitem{CSRL}
S.~Meyn.
\newblock {\em Control Systems and Reinforcement Learning}.
\newblock {Cambridge University Press}, Cambridge, 2022.

\bibitem{miller1990cmas}
W.~T. Miller, F.~H. Glanz, and L.~G. Kraft.
\newblock {CMAS}: An associative neural network alternative to backpropagation.
\newblock {\em Proceedings of the IEEE}, 78(10):1561--1567, 1990.

\bibitem{rahimi2007random}
A.~Rahimi and B.~Recht.
\newblock Random features for large-scale kernel machines.
\newblock {\em Advances in neural information processing systems}, 20, 2007.

\bibitem{rahimi2008weighted}
A.~Rahimi and B.~Recht.
\newblock Weighted sums of random kitchen sinks: Replacing minimization with
  randomization in learning.
\newblock {\em Advances in neural information processing systems}, 21, 2008.

\bibitem{schsei85}
P.~J. Schweitzer and A.~Seidmann.
\newblock Generalized polynomial approximations in {Markovian} decision
  processes.
\newblock {\em Journal of mathematical analysis and applications},
  110(2):568--582, 1985.

\bibitem{sutbar18}
R.~Sutton and A.~Barto.
\newblock {\em Reinforcement Learning: An Introduction}.
\newblock MIT Press, Cambridge, MA, 2nd edition, 2018.

\bibitem{sutszemae08}
R.~S. Sutton, C.~Szepesv\'{a}ri, and H.~R. Maei.
\newblock A convergent {O($n$)} algorithm for off-policy temporal-difference
  learning with linear function approximation.
\newblock In {\em Proc. of the Intl. Conference on Neural Information
  Processing Systems}, pages 1609--1616, Red Hook, NY, 2008.

\bibitem{tsivan96}
J.~N. Tsitsiklis and B.~Van~Roy.
\newblock Feature-based methods for large scale dynamic programming.
\newblock {\em Mach. Learn.}, 22(1-3):59--94, 1996.

\bibitem{wanboy09}
Y.~Wang and S.~Boyd.
\newblock Performance bounds for linear stochastic control.
\newblock {\em Systems Control Lett.}, 58(3):178--182, 2009.

\end{thebibliography}
